\documentclass[american]{article}
\usepackage[T1]{fontenc}
\usepackage[utf8]{inputenc}
\usepackage{babel}
\usepackage{algorithm2e}
\usepackage{amsmath}
\usepackage{amsthm}
\usepackage{amssymb}
\usepackage{geometry}
\usepackage{wasysym}
\usepackage[all]{xy}
\usepackage[]
 {hyperref}

\makeatletter
\numberwithin{figure}{section}
\numberwithin{table}{section}

\usepackage{tikz-cd}
\newsavebox{\diagOiiiCii}
\newsavebox{\diagOiiCiii}
\newsavebox{\diagOINii}
\newsavebox{\diagOiCiiFromI}
\newsavebox{\diagOICii}
\newsavebox{\diagOiCii}
\newsavebox{\diagI}
\newsavebox{\diagOiiiCiiB}
\newsavebox{\diagOiCiiB}

\newcommand{\begintikzsep}{\begin{tikzcd}[column sep = small, row sep = tiny]}

\makeatother

\begin{document}
\global\long\def\bgrp#1#2#3#4{\save[#1].[#2]!C#3*#4\frm{\{}*\frm{\}}\restore}%

\global\long\def\llp{\mathbin{\rightthreetimes}}%
\global\long\def\hom{\mathbin{\longrightarrow}}%
\global\long\def\sdiff{\mathbin{\backslash}}%

\xymatrixrowsep{5mm}
\xymatrixcolsep{5mm}
\title{Normal Spaces via Urysohn's Lemma as a Lifting Property}
\author{Robert Maxton}
\maketitle
\begin{abstract}
We present a translation of Urysohn's description of normal spaces
(as those where disjoint closed subsets are separated by a continuous
function) into the language of lifting properties in $\mathbf{Top}$,
correcting a frequently-cited erroneous translation. We also present
a translation of the definition of hereditarily normal spaces as those
in which every open subspace is normal, which we consider simpler
and more directly derived from the classical definition.
\end{abstract}

\section{Introduction}

In \cite{lift}, Misha Gavrilovich provides equivalent restatements
of many properties of topological spaces and continuous maps between
them in terms of categorical \emph{lifting properties} -- families
of commutative squares of the form
\[
\xymatrix{A\ar[rr]^{\psi}\ar[dd]_{f} &  & C\ar[dd]^{g}\\
\\B\ar[rr]_{\varphi}\ar[rruu]^{\exists\lambda} &  & D
}
\]
where we say that $f$ has the left lifting property against $g$,
or $g$ has the right lifting property against $f$, or $f\llp g$,
iff for every possible choice of $\psi,\varphi$ there exists some
morphism $\lambda:B\to C$ that makes the above diagram commute. In
particular, lifting properties are provided for most of the standard
separation axioms.

The author also introduces a convenient and concise notation for describing
continuous maps between specialization topologies: given one topology
presented as a graph of its specialization preorder, we simply repeat
the vertex labels to indicate where each point of the domain gets
sent to, using $=$ edges as necessary when two points are sent to
the same point in the codomain and using new vertex labels only for
points in the codomain with empty preimage. Continuity of the map
then follows by monotonicity of the map, i.e. directed edges are sent
to edges in the same direction. As an example, suppose that we define
the five point topology $\tau=\xymatrix{L\ar[r] & 0 & M\ar[l]\ar[r] & 1 & R\ar[l]}
$ and the topology of its dual order $\tau^{\dagger}:=\xymatrix{L^{\dagger} & 0^{\dagger}\ar[l]\ar[r] & M^{\dagger} & 1^{\dagger}\ar[l]\ar[r] & R^{\dagger}}
$; then the following:

\begin{lrbox}{\diagOiiiCii}%
  \begintikzsep
    L\arrow[rd] & & M\arrow[ld]\arrow[rd] & & R\arrow[ld] \\
      & 0 & & 1 &
  \end{tikzcd}
\end{lrbox}
\begin{lrbox}{\diagOiiCiii}%
\begintikzsep
	& L\arrow[ld]\arrow[rd] &&&& R\arrow[ld]\arrow[rd] & \\
	L^\prime && 0 \arrow[r, equal] & M\arrow[r, equal] & 1 && R^\prime
\end{tikzcd}
\end{lrbox}
\[
\xymatrix{
\left\{ \usebox\diagOiiiCii\right\} \ar@{|->}[r] & \left\{ \usebox\diagOiiCiii\right\}
}\]

represents the function
\begin{align*}
f: & \tau\to\tau^{\dagger}\\
 & x\mapsto\begin{cases}
0^{\dagger} & x=L\\
1^{\dagger} & x=R\\
M^{\dagger} & \text{otherwise}
\end{cases}
\end{align*}
where continuity is immediately apparent from the diagram.

This paper focuses on two translations of separation axioms into lifting
properties: Urysohn's lemma characterizing normal spaces, and the
definition of hereditarily normal spaces. Urysohn's lemma states that
a space $X$ is normal iff every pair of disjoint closed sets are
separated by a continuous function; that is, for every $s,t$ disjoint
and closed, there exists a continuous function $f:X\to\mathbb{R}$
s.t. $f(s)=\{0\},f(t)=\{1\}$. Let us abuse notation slightly; though
the unit interval is not a specialization topology, we nevertheless
represent the inclusion of the unit interval into a specialization
preorder with $\xymatrix@1{0\ \ar@{|-|}[r] &\ 1}$, and define a modified
unit interval with a topology roughly given by

\begin{lrbox}{\diagOINii}%
  \begintikzsep
    0^\prime\arrow[r, leftrightarrow] & 0\arrow[r, mapsto, no head]\arrow[r, mapsfrom, no tail] & 1\arrow[r, leftrightarrow] & 1^\prime
  \end{tikzcd}
\end{lrbox}
\[\left\{\usebox\diagOINii\right\}.\]

That is, to the original closed interval $\left[0,1\right]$, we 'glue
on' an extra $0,1$ topologically indistinguishable from their originals.
\cite{lift} then claims\footnote{In fact, the author claims that $\varnothing\to X$ lifts against
$\cdots\mapsto\left\{ \xymatrix{0^{\prime}\ar@{=}[r] & 0\ar[r] & [0,1] & 1\ar[l]\ar@{=}[r] & 1^{\prime}}
\right\} $. Reversing the arrows and assuming that the author intended to send
only the open interval $\left(0,1\right)$ to the central point makes
this continuous, but perhaps there is another intended reading.} that Urysohn's lemma is equivalent to the lifting property

\begin{lrbox}{\diagOiCiiFromI}%
  \begintikzsep
      & & (0, 1) \arrow[ld]\arrow[rd] \\
    0^\prime \arrow[r, equal] & 0 & & 1 \arrow[r, equal] & 1^\prime
  \end{tikzcd}
\end{lrbox}
\[\xymatrix{\varnothing\ar[r] & X\ar@{}[r]|(0.2)*+{\llp} & \left\{\usebox\diagOINii\right\}\ar@{|->}[r] & \left\{\usebox\diagOiCiiFromI\right\}}.\]However, this does not quite work. We now provide a slightly modified
translation; in proving that our modification faithfully reproduces
Urysohn's lemma, we will also clarify why the original does not.

\section{A Modification}

Let us define the 'interval with doubled endpoints' $I_{0}^{1}$,
slightly modified from the original paper: we have the closed unit
interval $[0,1]$, onto which we have again glued a $0^{\prime},1^{\prime}$,
but here we say that $0\leadsto0^{\prime}$ and $1\leadsto1^{\prime}$
only, rather than making them indistinguishable. We may then use Urysohn's
lemma to define a lifting property for normal spaces: a space $X$
is normal iff the following square has a lift for all $\chi$:

\begin{lrbox}{\diagOICii}%
  \begintikzsep
      & 0 \arrow[r, mapsto, no head]\arrow[r, mapsfrom, no tail]\arrow[ld] & 1\arrow[rd] & \\
    0^\prime & & & 1^\prime
  \end{tikzcd}
\end{lrbox}
\begin{lrbox}{\diagOiCii}%
  \begintikzsep
      & 0 \arrow[r, equal]\arrow[ld] & 1\arrow[rd] & \\
    0^\prime & & & 1^\prime
  \end{tikzcd}
\end{lrbox}
\begin{lrbox}{\diagI}%
  \begintikzsep
    0^\prime \arrow[r, equal] & 0 \arrow[r, mapsto, no head]\arrow[r, mapsfrom, no tail] & 1\arrow[r, equal] & 1^\prime
  \end{tikzcd}
\end{lrbox}
\[
\xymatrix{
\varnothing\ar[rr]\ar[d] & & \left\{ \usebox\diagOICii\right\} \ar@{|->}[d]^{\pi}\ar@{|->}[r]^{\iota} & \left\{ \usebox\diagI\right\} \\X\ar[rr]^(0.3)\chi\ar@{-->}[rru]^{\exists\lambda} & & \left\{ \usebox\diagOiCii\right\} 
}\]
\begin{proof}
We may read off a given choice of two closed sets in $X$ as the fibers
$s\equiv\chi^{-1}\left(0^{\prime}\right),t\equiv\chi^{-1}\left(1^{\prime}\right)$
of the closed points $0^{\prime},1^{\prime}$ under the bottom arrow
$\chi$; a lift $\lambda$ through $\pi$ then gives us a continuous
map to the modified interval $I_{0}^{1}$, with our two closed sets
now the preimages $\lambda^{-1}\left(0^{\prime}\right),\lambda^{-1}\left(1^{\prime}\right)$
in $I_{0}^{1}$ by commutativity of the diagram. By composing with
$\iota=\left\{ \xymatrix{0^{\prime} & 0\ar[l]\ar@{|-|}[r] & 1\ar[r] & 1^{\prime}}
\right\} \mapsto\left\{ \xymatrix{0^{\prime}=0\ar@{|-|}[r] & 1=1^{\prime}}
\right\} $, we have a continuous map to the usual interval $I$, where $\lambda^{-1}\left(0^{\prime}\right)\subseteq\left(\iota\circ\lambda\right)^{-1}\left(0\right)$
and $\lambda^{-1}\left(1^{\prime}\right)\subseteq\left(\iota\circ\lambda\right)^{-1}\left(1\right)$,
or in other words that $s,t$ are separated by the continuous map
$\iota\circ\lambda$; by Urysohn's lemma $X$ is normal. 

Conversely, suppose $X$ is normal, and therefore for any pair of
disjoint closed sets $s,t$ we have a separating continuous map $\lambda^{\prime}:X\to I$.
Then we may define
\begin{align*}
\lambda\left(x\right)= & \begin{cases}
0^{\prime} & x\in s\\
0 & \lambda^{\prime}\left(x\right)=0\wedge x\notin s\\
1^{\prime} & x\in t\\
1 & \lambda^{\prime}\left(1\right)=1\wedge x\notin t\\
\lambda^{\prime}\left(x\right) & \text{otherwise}
\end{cases}
\end{align*}
whence $\lambda^{\prime}=\iota\circ\lambda$ and choosing $\chi=\pi\circ\lambda$
makes the above lifting square commute.
\end{proof}
Note that for this construction to work, $I_{0}^{1}$ needed a morphism
$\pi$ to $\left\{\xymatrix@1{L & 0\ar[l]\ar[r] & R}\right\}$ and
a morphism $\iota$ to $\left[0,1\right]$ such that the fiber $\pi^{-1}\left(L\right)$
is a \emph{strict} subset $\subsetneq\iota^{-1}\left(0\right)$, and
similarly $\pi^{-1}\left(R\right)\subsetneq\iota^{-1}\left(1\right)$.
If, as in \cite{lift}, we had defined $I_{0}^{1}$ with $0^{\prime}\leftrightarrow0$
and $1\leftrightarrow1^{\prime}$, then the image of $\left[0,1\right]$
would not be open in $I_{0}^{1}$, we would be obliged to send only
the open interval $\left(0,1\right)$ to the open central point in
$\left\{\xymatrix@1{L & 0\ar[l]\ar[r] & R}\right\}$, and we would
have $\pi^{-1}\left(L\right)=\iota^{-1}\left(0\right)$ and similarly
for $R,1$. In this case we instead reproduce the translation of perfectly
normal spaces, where disjoint closed sets are precisely separated
by a continuous function.

\section{A Simplification}

As an addendum, we also note that the lifting property translating
hereditarily normal spaces may be simplified. Recall that a hereditarily
normal space may be defined as one where every open subspace is itself
normal. Then we may transform the lifting property for normal spaces
into one for hereditarily normal spaces mechanically, by 'hoisting'
the defining morphism atop a new, closed point that becomes the target
of the complement of each considered open set. That is, the lifting
square given by

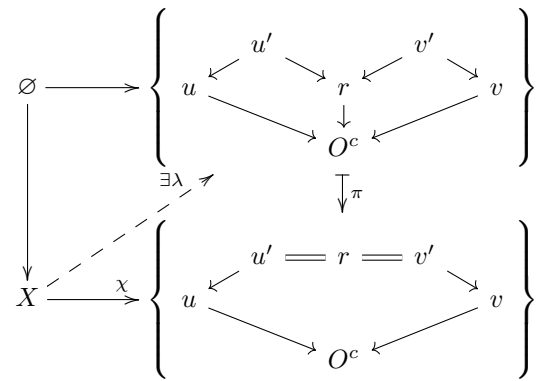
\begin{lrbox}{\diagOiiiCiiB}%
  \begintikzsep
      & u^\prime \arrow[ld]\arrow[rd] &  & v^\prime \arrow[ld]\arrow[rd] & \\
    u \arrow[rrdd] &  & r\arrow[dd] & & v \arrow[lldd] \\ \\
	  &  & O^c
  \end{tikzcd}
\end{lrbox}
\begin{lrbox}{\diagOiCiiB}%
  \begintikzsep
      & u^\prime \arrow[r,equal] \arrow[ld] & r \arrow[r, equal] & v^\prime \arrow[rd] & \\
    u \arrow[rrdd] &  & & & v \arrow[lldd] \\ \\
	  &  & O^c
  \end{tikzcd}
\end{lrbox}
\[
\xymatrix{
\varnothing\ar[rr]\ar[d] & & \left\{ \usebox\diagOiiiCiiB\right\} \ar@{|->}[d]^{\pi} \\X\ar[rr]^(0.3)\chi\ar@{-->}[rru]^{\exists\lambda} & & \left\{ \usebox\diagOiCiiB\right\} 
}\]

directly translates the definition.

\section{Conclusion}

We have here provided a pair of translations of standard topological
separation axioms into categorical language. A formalization in Lean,
making use of Mathlib, can be found at \href{https://github.com/robertmaxton42/mathlib4/tree/separation}{https://github.com/robertmaxton42/mathlib4/tree/separation}.

\end{document}